\documentclass{article}
\usepackage{amssymb}
\usepackage{amsfonts}
\usepackage{amsmath}
\usepackage{doublespace}
\usepackage{acronym}

\newtheorem{theorem}{Theorem}

\newtheorem{lemma}[theorem]{Lemma}

\newtheorem{proposition}[theorem]{Proposition}
\newtheorem{remark}[theorem]{Remark}

\newenvironment{proof}[1][Proof]{\noindent\textbf{#1.} }{\ \rule{0.5em}{0.5em}}
\input{tcilatex}

\begin{document}

\title{Ratio vectors of fourth degree polynomials}
\author{Alan Horwitz \\
%EndAName
Penn State University\\
25 Yearsley Mill Rd.\\
Media, PA 19063\\
alh4@psu.edu}
\maketitle
\date{}

\begin{abstract}
Let $p(x)$ be a polynomial of degree $4$ with four distinct real roots $%
r_{1}<r_{2}<r_{3}<r_{4}$. Let $x_{1}<x_{2}<x_{3}$ be the critical points of $%
p$, and define the ratios $\sigma _{k}=\dfrac{x_{k}-r_{k}}{r_{k+1}-r_{k}}%
,k=1,2,3$. For notational convenience, let $\sigma _{1}=u$, $\sigma _{2}=v$,
and $\sigma _{3}=w$. $(u,v,w)$ is called the ratio vector of $p.$ We prove
necessary and sufficient conditions for $(u,v,w)$ to be a ratio vector of a
polynomial of degree $4$ with all real roots. Most of the necessary
conditions were proven in (\cite{h}). The main results of this paper involve
using the theory of Groebner bases to prove that those conditions are also
sufficient.

\end{abstract}

{}{\LARGE Introduction}

\qquad Let $p(x)$ be a polynomial of degree $n\geq 2$ with $n$ distinct real
roots $r_{1}<r_{2}<\cdots <r_{n}$. Let $x_{1}<x_{2}<\cdots <x_{n-1}$ be the
critical points of $p$, and define the ratios

\begin{equation}
\sigma _{k}=\dfrac{x_{k}-r_{k}}{r_{k+1}-r_{k}},k=1,2,...,n-1  \label{ratio}
\end{equation}

$(\sigma _{1},...,\sigma _{n-1})$ is called the ratio \textbf{vector} of $p$%
, and $\sigma _{k}$ is called the $k$th ratio.The following inequality(\ref%
{ineq}) was derived in (\cite{a}) and cited in the author's paper (\cite{h}%
). I have since discovered that it was actually derived ealier by Peyser in (%
\cite{p}).

$\medskip $%
\begin{equation}
\frac{1}{n-k+1}\,<\,\sigma _{k}<\frac{k}{k+1}\,\,  \label{ineq}
\end{equation}

Solving $\sigma _{k}=\dfrac{x_{k}-r_{k}}{r_{k+1}-r_{k}}$ for $x_{k}$ yields 
\begin{equation}
x_{k}=(r_{k+1}-r_{k})\sigma _{k}+r_{k},\;k=1,2,...,n-1  \label{xksk}
\end{equation}

Andrews also defined the sets $X_{n}=\prod\limits_{k=1}^{n-1}\left( \dfrac{1%
}{n-k+1},\dfrac{k}{k+1}\right) $ and $Y_{n}=$ the set of elements in $%
R^{n-1} $ which are the ratio vectors of polynomials with $n$ distinct real
zeroes. Our main results in this paper concern the case $n=4$. For
notational convenience, let $\sigma _{1}=u$, $\sigma _{2}=v$, and $\sigma
_{3}=w$. Then (\ref{ineq}) becomes

\begin{equation}
\text{\ }X_{4}=\left\{ (u,v,w):\frac{1}{4}<u<\frac{1}{2}\text{, }\frac{1}{3}%
<v<\frac{2}{3}\text{,}\frac{1}{2}<w<\frac{3}{4}\right\}  \label{x4}
\end{equation}%
In (\cite{h}) we showed that $Y_{4}$ lies on the zero set of a polynomial of
degree $9$, $Q(u,v,w)=(1-4v+4uv)R(u,v,w)$, where $R(u,v,w)$ is defined below(%
\ref{R}), which is a \textit{necessary} condition for $(u,v,w)$ to be a
ratio vector. The polynomial $Q$ in (\cite{h}) was found by directly solving
the system of equations (\ref{1})--(\ref{3}) below. This is not difficult to
do by solving for $r$ in terms of $s$ in the first two equations and then
substituting into the third equation(see \cite{h}).We show below that $%
1-4v+4uv<0$ for $(u,v,w)\in X_{4}$, so that $Y_{4}$ actually lies on the
zero set, $Z$, of $R$. Thus a \textbf{necessary} condition for $(u,v,w)$ to
be a ratio vector is $(u,v,w)\in Z\cap X_{4}$. One can also show that a
ratio vector must satisfy certain more restrictive inequalities than those
which define $X_{4}$(see (\ref{5})--(\ref{7})). This turns out to be
equivalent to the condition that $(u,v,w)\in Z\cap \left\{ Z_{1}\cup
Z_{2}\cup Z_{3}\right\} .$ Our main result(see Theorem \ref{main}) is that $%
(u,v,w)\in Z\cap \left\{ Z_{1}\cup Z_{2}\cup Z_{3}\right\} $ is also\textit{%
\ }\textbf{sufficient} for $(u,v,w)$ to be a ratio vector. To obtain \textit{%
sufficient} conditions, however, one needs more sophisticated methods, like
the theory of \textbf{Gr\"{o}ebner bases}, which we use throughout the paper.

{\LARGE Preliminary Material}

We shall need a system of algebraic equations which relate the real roots, $%
r_{k},$ and the ratios, $\sigma _{k}$. One can rewrite $p(x)=(x-r_{1})\cdots
(x-r_{n})$ using the elementary symmetric functions $E_{j}=\dfrac{e_{j}}{%
\binom{n}{j}}$, $e_{j}\equiv e_{j}(r_{1},...,r_{n})=j$th elementary
symmetric function of the $r_{j}$, $j=1,2,...,n$, starting with $%
e_{1}(r_{1},...,r_{n})=r_{1}+\cdots r_{n}$, etc. 
\begin{equation}
p(x)=\ \dsum\limits_{j=0}^{n-1}(-1)^{j}\binom{n}{j}E_{j}x^{n-j}+x^{n}
\label{4}
\end{equation}%
$\,$Let $x_{1}<x_{2}<\cdots <x_{n-1}$ denote the critical points of $p$. By
equating coefficients it is easy to show(see \cite{r}) that 
\begin{equation}
E_{j}\,(r_{1},...,r_{n})=E_{j}\,(x_{1},...,x_{n-1}),\;j=1,2,...,n-1
\label{symm}
\end{equation}

We shall also make use of

\begin{lemma}
\label{cp}Suppose that (\ref{symm}) holds for distinct numbers $%
r_{1}<r_{2}<\cdots <r_{n}$ and $x_{1}<x_{2}<\cdots <x_{n-1}$. Let $%
p(x)=(x-r_{1})(x-r_{2})\cdots (x-r_{n})$. Then the $x_{j}$ must be the
critical points of $p$.
\end{lemma}

%TCIMACRO{\TeXButton{Proof}{\proof}}%
%BeginExpansion
\proof%
%EndExpansion
Since $\binom{n}{j}(n-j)=\dfrac{n!}{j!(n-j-1)!}=n\binom{n-1}{j},$ by (\ref{4}%
) $p^{\prime }(x)=\dsum\limits_{j=1}^{n-1}(-1)^{j}\binom{n}{j}%
(n-j)E_{j}(r_{1},...,r_{n})x^{n-j-1}+nx^{n-1}=n\dsum%
\limits_{j=1}^{n-1}(-1)^{j}\binom{n-1}{j}%
E_{j}(r_{1},...,r_{n})x^{n-j-1}+nx^{n-1}=$

$n\dsum\limits_{j=1}^{n-1}(-1)^{j}\binom{n-1}{j}%
E_{j}(x_{1},...,x_{n-1})x^{n-j-1}+nx^{n-1}=n(x-x_{1})(x-x_{2})\cdots
(x-x_{n-1})$, which implies that the $x_{j}$ must be the critical points of $%
p$.

We shall now consider the case $n=4$ for the rest of the paper. If $p$ is a
fourth degree polynomial with four distinct real zeroes, then $p(x+c)$ and $%
p(cx)(c\neq 0)$ have the same ratio vectors as $p$. Thus we may assume that
the zeros of $p$ are $r_{1}=-1<$ $r_{2}=0<r_{3}<r_{4}$. Also, for notational
convenience, let $r_{3}=r$, $r_{4}=s$. Then (\ref{symm}) becomes 
\begin{eqnarray}
4x_{1}x_{2}x_{3}+rs &=&\allowbreak 0  \notag \\
4(x_{1}+x_{2}+x_{3})-3(-1+r+s) &=&\allowbreak 0  \label{symm4} \\
2(x_{1}x_{2}+x_{1}x_{3}+x_{2}x_{3})-(-r-s+rs) &=&\allowbreak 0  \notag
\end{eqnarray}

Our system of algebraic equations will come from substituting for the $x_{k}$
in (\ref{symm4}) using (\ref{xksk}). Using $x_{1}=u-1$, $x_{2}=rv$, $%
x_{3}=(s-r)w+r$, we have, after some cancellation, the equivalent system of
equations 
\begin{equation}
\left( 4(1-u)vw-1\right) s+4(1-u)v(1-w)r=0  \label{1}
\end{equation}%
\begin{equation}
\left( -1+4u\right) +\left( 1+4v-4w\right) r+\left( 4w-3\right) s=0
\label{2}
\end{equation}%
\begin{equation}
\left( 2u(v-w+1)-2v+2w-1\right) r+(2uw-2w+1)s\,-2v(w-1)r^{2}-(1-2vw)rs=0
\label{3}
\end{equation}

{\LARGE Main Results}

First we state some inequalities for ratio vectors of fourth degree
polynomials which are of interest in their own right, and which will also be
used to prove our main result below.

\begin{lemma}
\label{L1}Suppose that $(u,v,w)$ is a ratio vector of a fourth degree
polynomial with four distinct real zeroes. Then 
\begin{equation}
\dfrac{1}{4(1-u)}<v<\dfrac{1}{4(1-u)w}  \label{5}
\end{equation}

\begin{equation}
\dfrac{1}{4(1-v)}<w<\dfrac{1}{4(1-u)(1-v)}  \label{6}
\end{equation}

and 
\begin{equation}
w<\dfrac{1}{2(1-u)}  \label{7}
\end{equation}
\end{lemma}

\begin{proof}
As noted above, we may assume that the zeros of $p$ are $-1<0<r<s$. Since $%
r<s$ and $(u,v,w)\in X_{4}$, by (\ref{1}), $0=\left( 4(1-u)vw-1\right)
s+4(1-u)v(1-w)r<\left( (4(1-u)vw-1)s+4(1-u)v(1-w)\right) s=\allowbreak
(-1+4v-4vu)s$. Thus 
\begin{equation}
-1+4v-4vu>0  \label{5a}
\end{equation}%
, which is equivalent to $v>\dfrac{1}{4(1-u)}$ since $(u,v,w)\in
X_{4}\Rightarrow 1-u>0$. Again, by (\ref{1}), $\dfrac{s}{r}=\dfrac{%
4(1-u)v(1-w)}{-4(1-u)vw+1}=\dfrac{\allowbreak 4(1-u)v(1-w)}{-4vw+4vwu+1}$.
Since $\dfrac{s}{r}>0$ and $\allowbreak 4(1-u)v(1-w)>0$%
\begin{equation}
-4vw+4vwu+1>0  \label{5b}
\end{equation}%
, which is equivalent to $v<\dfrac{1}{4(1-u)w}$. That proves (\ref{5}). (\ref%
{6}) follows in a similar fashion by assuming that $%
r_{1}<r_{2}<r_{3}=0<r_{4} $, which yields a system analogous to (\ref{1})--(%
\ref{3}). To prove (\ref{7}), note that the upper bounds in (\ref{5}) and (%
\ref{6}) imply that $4(1-u)vw-1<0$ and $4(1-u)(1-v)w-1<0$. Thus $\dfrac{1}{2}%
\left( 4(1-u)vw-1+4(1-u)(1-v)w-1\right) <0\Rightarrow $

\begin{equation}
-1+2w-2wu<0  \label{7a}
\end{equation}%
, which is equivalent to $w<\dfrac{1}{2(1-u)}$.
\end{proof}

\begin{remark}
We note here that the lower bounds in (\ref{5}) and in (\ref{6}) can be used
to prove the monotonicity of the ratios. That is, to prove that $u<v$(see (%
\cite{h})) and that $v<w.$
\end{remark}

Define the polynomial in three variables

\begin{gather}
R(u,v,w)=32u^{2}v^{3}w^{2}-48u^{2}v^{2}w^{2}-64uv^{3}w^{2}+24u^{2}v^{2}w+16u^{2}vw^{2}+120uv^{2}w^{2}
\notag \\
+32v^{3}w^{2}-16u^{2}vw-48uv^{2}w-64uvw^{2}-72v^{2}w^{2}+52uvw+8uw^{2}
\label{R} \\
+24v^{2}w+48vw^{2}-6uv-10uw-30vw-8w^{2}+2u+4v+6w-1  \notag
\end{gather}

along with its zero set

\begin{equation*}
Z=\left\{ (u,v,w):R(u,v,w)=0\right\}
\end{equation*}

Also define the following subsets of $X_{4}$: 
\begin{equation*}
Z_{1}=\left\{ (u,v,w):\dfrac{1}{4}<u\leq \dfrac{1}{3},\dfrac{1}{4(1-u)}<v<%
\dfrac{1}{2},\dfrac{1}{2}<w<\dfrac{1}{4(1-u)(1-v)}\right\}
\end{equation*}%
\begin{equation*}
Z_{2}=\left\{ (u,v,w):\dfrac{1}{4}<u\leq \dfrac{1}{3},\dfrac{1}{2}<w<\dfrac{1%
}{2(1-u)},\dfrac{1}{2}\leq v<\dfrac{1}{4(1-u)w}\right\} \cap \left\{
(u,v,w):v<\dfrac{2}{3}\right\}
\end{equation*}

\begin{equation*}
Z_{3}=\left\{ (u,v,w):\dfrac{1}{3}<u<\dfrac{1}{2},\dfrac{1}{2}<w<\dfrac{3}{4}%
,\dfrac{1}{4(1-u)}<v<\dfrac{1}{4(1-u)w}\right\} \cap \left\{ (u,v,w):v<%
\dfrac{2}{3}\right\}
\end{equation*}

We now give necessary and sufficient conditions for $(u,v,w)$ to be a ratio
vector.

\begin{theorem}
\label{main}$(u,v,w)$ is a ratio vector of a fourth degree polynomial with
four distinct real zeroes if and only if $(u,v,w)\in Z\cap \left\{ Z_{1}\cup
Z_{2}\cup Z_{3}\right\} $.
\end{theorem}

Thus the set of ratio vectors of fourth degree polynomials is precisely
equal to $Z\cap \left\{ Z_{1}\cup Z_{2}\cup Z_{3}\right\} $.

Theorem \ref{main} will follow directly from the following two propositions.
Most of Proposition \ref{necessary} was proved in (\cite{h}).

\begin{proposition}
\label{necessary}Suppose that $(u,v,w)$ is a ratio vector of a fourth degree
polynomial with four real distinct zeros. Then $(u,v,w)\in Z\cap \left\{
Z_{1}\cup Z_{2}\cup Z_{3}\right\} $.
\end{proposition}

\begin{proof}
As noted earlier, we showed in (\cite{h}, Theorem 2) that $Y_{4}$ is
contained in the zero set of the polynomial $Q(u,v,w)=(1-4v+4uv)R(u,v,w)$.
By (\ref{5a}), $R(u,v,w)=0$, and thus $(u,v,w)\in Z$. $R(u,v,w)=0$ will also
follow from the Gr\"{o}ebner basis we use below in the proof of Proposition %
\ref{sufficient}. $(u,v,w)\in Z_{1}\cup Z_{2}\cup Z_{3}$ follows immediately
from (\ref{x4}) and Lemma \ref{L1}.
\end{proof}

\textbf{Example: }Let $p(x)=\allowbreak
128x^{4}-752x^{3}+1636x^{2}-1558x+546 $, so that the roots are $r_{1}=1$, $%
r_{2}=\dfrac{3}{2}$, $r_{3}=\dfrac{13}{8}$, $r_{4}=\dfrac{7}{4}$, and the
critical points are $x_{1}\approx 1.1506$, $x_{2}\approx 1.\allowbreak 5560$%
, $x_{3}\approx 1.\allowbreak 6996$. The ratios are $u=\dfrac{x_{1}-r_{1}}{%
r_{2}-r_{1}}\approx \ .301\,3$, $v=\dfrac{x_{2}-r_{2}}{r_{3}-r_{2}}\approx \
.448\,1$, $w=\dfrac{x_{3}-r_{3}}{r_{4}-r_{3}}\approx \ .596\,8$. Then $u<%
\dfrac{1}{3}$, $\dfrac{1}{4(1-u)}\approx \ .357\,8<v<\dfrac{1}{2}$, and $%
\dfrac{1}{2}<w<\dfrac{1}{4(1-u)(1-v)}\approx \ .648\,3$, so that $(u,v,w)\in
Z_{1}$.

We prove a simple lemma about $Z_{1}\cup Z_{2}\cup Z_{3}$.

\begin{lemma}
\label{xz}If $(u,v,w)\in \left\{ Z_{1}\cup Z_{2}\cup Z_{3}\right\} $, then $%
(u,v,w)\in X_{4}$ and (\ref{5}) holds, or equivalently, (\ref{5a}) and (\ref%
{5b}).
\end{lemma}

\begin{proof}
Case 1: $(u,v,w)\in $ $Z_{1}$. Then $v>\dfrac{1}{4(1-u)}>\dfrac{1}{3}$ and $%
w<\dfrac{1}{4(1-u)(1-v)}<\dfrac{3}{4}$, which implies that $(u,v,w)\in X_{4}$%
. Also, $\dfrac{1}{4(1-u)w}>1-v>v$, so that (\ref{5}) holds.

Case 2: $(u,v,w)\in $ $Z_{2}$. Then $w<\dfrac{1}{2(1-u)}<\dfrac{3}{4}$,
which implies that $(u,v,w)\in X_{4}$. Also, $\dfrac{1}{4(1-u)}<\dfrac{3}{8}%
<v$, so that (\ref{5}) holds.

Case 3: $(u,v,w)\in $ $Z_{3}$. Then $v>\dfrac{1}{4(1-u)}>\dfrac{3}{8}>\dfrac{%
1}{3}$ which implies that $(u,v,w)\in X_{4}$. It is obvious that (\ref{5})
holds.
\end{proof}

Before proving the sufficiency part of Theorem \ref{main} in the form of
Proposition \ref{sufficient} below, we define the following polynomial,
which will be important in our proof. 
\begin{equation}
k(u,v,w)=8u^{2}v^{2}w-16uv^{2}w+8uvw+8v^{2}w-2uv-2u\allowbreak
w-8vw+2u+2v+2w-1  \label{k}
\end{equation}%
We now prove the following key lemma about $k$.

\begin{lemma}
\label{dk}If $(u,v,w)\in Z\cap \left\{ Z_{1}\cup Z_{2}\cup Z_{3}\right\} $,
then $k(u,v,w)>0$.
\end{lemma}

\begin{proof}
Write $k(u,v,w)=\allowbreak \allowbreak -2\left( 1-u\right) \left(
4v^{2}u-4v^{2}+4v-1\right) w-2uv+\allowbreak 2u+2v-1$, and let $f(v)=\dfrac{%
4v^{2}-4v+1}{4v^{2}}$. Then, for $\dfrac{1}{3}<v<\dfrac{2}{3},$ $f\,^{\prime
\prime }(v)>0$. Since $f\left( \dfrac{1}{3}\right) =\allowbreak \dfrac{1}{4}$
and $f\left( \dfrac{2}{3}\right) =\allowbreak \dfrac{1}{16}$, $f(v)\leq 
\dfrac{1}{4}$ for $\dfrac{1}{3}<v<\dfrac{2}{3}.$ Thus $u>\dfrac{1}{4}\geq 
\dfrac{1+4v^{2}-4v}{4v^{2}}$, which implies that $4v^{2}u-4v^{2}+4v-1>0$.
That proves that $k(u,v,w)$ is a \textbf{decreasing} function of $w$ for any 
$\dfrac{1}{4}<u<\dfrac{1}{2}$, $\dfrac{1}{3}<v<\dfrac{2}{3}$. For $%
(u,v,w)\in Z_{1}$, the largest value of $w$ is $\dfrac{1}{4(1-u)(1-v)}$. Now 
$k\left( u,v,\dfrac{1}{4(1-u)(1-v)}\right) =\allowbreak \dfrac{1}{2}\left(
1-2v\right) \dfrac{4u-1}{1-v}>0$ since $v<\dfrac{1}{2}$. Thus $k(u,v,w)>0$
on all of $Z_{1}$. For $(u,v,w)\in Z_{2}$, the largest value of $w$ is $%
\dfrac{1}{2(1-u)}$, so we consider $k\left( u,v,\dfrac{1}{2(1-u)}\right)
=\allowbreak 2\left( 2v-1\right) \left( v-uv-u\right) $. $u\leq \dfrac{1}{3}%
\Rightarrow \dfrac{u}{1-u}\leq \dfrac{1}{2}\leq v\Rightarrow v-uv-u\geq
0\Rightarrow k\left( u,v,\dfrac{1}{2(1-u)}\right) \geq 0$. Again, it follows
that $k(u,v,w)\geq 0$ on all of $Z_{2}$. Finally, suppose that $(u,v,w)\in
Z_{3}$. Then $2k\left( u,v,\dfrac{3}{4}\right) =H(u,v)\allowbreak $, where $%
\allowbreak H(u,v)=12u^{2}v^{2}-24uv^{2}+8uv+12v^{2}+u-8v+1$. Now $%
H(1/3,v)=\allowbreak \dfrac{4}{3}\left( 2v-1\right) ^{2}\geq 0$, $%
H(1/2,v)=\allowbreak 3v^{2}+\dfrac{3}{2}-4v>0$, $H(u,1/3)=\dfrac{1}{3}\left(
u+1\right) \left( 4u-1\right) >0$, and $H(u,2/3)=\allowbreak \dfrac{16}{3}%
u^{2}-\dfrac{13}{3}u+1>0$. In addition, setting $\dfrac{\partial H(u,v)}{%
\partial u}=\dfrac{\partial H(u,v)}{\partial v}=0$ yields only one solution, 
$u=1,v=-\dfrac{1}{8}$. Thus $H(u,v)\geq 0$ for $\dfrac{1}{3}<u<\dfrac{1}{2},%
\dfrac{1}{3}<v<\dfrac{2}{3}$, which implies that $k(u,v,w)\geq 0$ on all of $%
Z_{3}$. Thus $k(u,v,w)\geq 0$ on $Z_{1}\cup Z_{2}\cup Z_{3}$. Now let $%
d(u,v,w)=-16uv^{2}w+24uvw+16v^{2}w-12uv-24vw+\allowbreak 8v+4w-1$, and note
the following identity relating $d,k,$ and $R$: $d(u,v,w)k(u,v,w)=2v\left(
4u-1\right) \left( 1-w\right) \left( 1-2u\right) \left( -4vw+4vwu+1\right)
+\left( -1+4v-4vu\right) R(u,v,w)$

Then if $R(u,v,w)=0$, $k(u,v,w)\neq 0$ since $(u,v,w)\in Z_{1}\cup Z_{2}\cup
Z_{3}$(by (\ref{5b})). It follows immediately that $k(u,v,w)>0$ if $%
(u,v,w)\in Z\cap \left\{ Z_{1}\cup Z_{2}\cup Z_{3}\right\} $.
\end{proof}

\begin{remark}
With a little more effort, one could actually prove that $k(u,v,w)>0$ if $%
(u,v,w)\in Z_{1}\cup Z_{2}\cup Z_{3}$, but we did not require that result.
Note that $k(u,v,w)$ is \textbf{not} nonnegative on all of $X_{4}$. For
example, $k(9/32,4/9,7/10)=\allowbreak -\dfrac{127}{12\,960}$.
\end{remark}

We will now prove the sufficiency part of Theorem \ref{main} in the form of
the following proposition.

\begin{proposition}
\label{sufficient} Suppose that $(u,v,w)\in Z\cap \left\{ Z_{1}\cup
Z_{2}\cup Z_{3}\right\} $. Then there are unique real numbers $0<r<s$ such
that the polynomial $p(x)=(x+1)x(x-r)(x-s)$ has $(u,v,w)$ as a ratio vector.
Furthermore, $r=\dfrac{k(u,v,w)}{\allowbreak 2v\left( 1-2u\right) \left(
1-w\right) }$ and $s=\dfrac{\allowbreak 2(1-u)vk(u,v,w)}{\left(
-4vw+4vwu+1\right) v\left( 1-2u\right) }$.
\end{proposition}

\begin{proof}
Let $C[t_{1},...,t_{n}]$ denote the polynomials in $t_{1},...,t_{n}$ with
complex coefficients, and for any ideal $I\subseteq C[t_{1},...,t_{n}],$ let 
$V(I)=\left\{ (t_{1},...,t_{n}):f(t_{1},...,t_{n})=0\;\forall f\in I\right\} 
$. Our approach is to obtain as much information as one can by viewing $%
s,r,u,v,w$ as independent variables in (\ref{1})--(\ref{3}), even though in
reality they are not by (\ref{xksk}). Let $f,g,$ and $h$ denote the LHS of
equations (\ref{1})--(\ref{3}),$\;I=$ $\langle $ $f,g,h\rangle =$ the ideal
generated by $f,g,h$ in $C[s,r,u,v,w]$. Let $I_{1}$be the first elimination
ideal, $I\cap C[r,u,v,w],$ and let $I_{2}$ equal the second elimination
ideal, $I\cap C[u,v,w]$. We found a Gr\"{o}ebner basis for $I$, denoted by $%
LEX$, using Maple 7 with the \textbf{lexographic ordering} $s>r>w>v>u$. $LEX$
contains $10$ elements, which we denote by $LEX_{1},...,LEX_{10}$. Since $I=$
$\langle LEX_{1},...,LEX_{10}\rangle ,$ it follows that equations (\ref{1}%
)--(\ref{3}) and the system of equations $LEX_{1}=0,...,LEX_{10}=0$ have
exactly the same set of solutions. $LEX_{1}$ is the only element of $LEX$
which only depends on $u,v,$ and $w$\footnote{%
Different computer algebra systems may list the elements of a minimal
reduced Grobner basis in a different order, but they must contain the same
polynomials.}. Hence, by the Elimination Theorem(see (\cite{clo}, Theorem 2,
page 114), $LEX_{1}$ is a Gr\"{o}ebner basis for the second elimination
ideal, $I_{2}$. Since $R(u,v,w)=0$ by assumption and $LEX_{1}$ is a multiple
of $R$, $(u,v,w)$ is a \textit{partial} solution of (\ref{1})--(\ref{3}) in $%
V\left( I_{2}\right) $. Let $S=\{LEX_{1},...,LEX_{6}\}$, \textbf{none} of
which involve $s$. Then by the Elimination Theorem again, $S$ is a Gr\"{o}%
ebner basis for the first elimination ideal, $I_{1}$. In particular, $%
LEX_{2}=v(4v-3)(2u-1)(4u-1)(1-4v+4uv)r+c(u,v,w)$ for some polynomial $%
c(u,v,w)$. Now $v(4v-3)(2u-1)(4u-1)(1-4v+4uv)\neq 0$ if $(u,v,w)\in
Z_{1}\cup Z_{2}\cup Z_{3}$ by Lemma \ref{xz}. Thus \textbf{all} of \ the
coefficients of the highest powers of $r$ in each element of $S$ \ \textbf{%
cannot} vanish. By the Extension Theorem(see (\cite{clo}, Theorem 3, page
117), for each $(u,v,w)\in Z\cap \left\{ Z_{1}\cup Z_{2}\cup Z_{3}\right\} $%
, there is a complex number $r$ such that $(r,u,v,w)\in V(I_{1})$\footnote{%
We use the fact that $I_{2}$ is the first elimination ideal of $I_{1}$.}.
Now $LEX_{8}=(4w-3)s+$ $r+4rv-4rw+4u-1$. Since $4w-3\neq 0$, by the
Extension Theorem again, there is a complex number $s$ such that $%
(s,r,u,v,w)\in V(I)$. That is, $(s,r,u,v,w)$ is a solution of (\ref{1})--(%
\ref{3}). The next thing we need to show is that $0<r<s$, and that $r$ and $%
s $ are unique for each given $(u,v,w)$.We now find it more convenient to
use a Gr\"{o}ebner basis for $I$ using the \textbf{total degree ordering},
which we denote by \textbf{\ }$TDEG$. $TDEG$ has $7$ elements and again,
equations (\ref{1})--(\ref{3}) and the system of equations $%
TDEG_{1}=0,...,TDEG_{7}=0$ have exactly the same set of solutions. In
particular, $TDEG_{7}=(4u-1)A(r,v,u,w)$, where $A(r,v,u,w)=\allowbreak
\allowbreak -2v\left( 2u-1\right) \left( w-1\right) r+$

$8u^{2}v^{2}w-16uv^{2}w+8uvw+8v^{2}w-2uv-2u\allowbreak w-8vw+2u+2v+2w-1=$

$-2v\left( 2u-1\right) \left( w-1\right) r+k(u,v,w)$. $TDEG_{7}=0\Rightarrow
-2v\left( 2u-1\right) \left( w-1\right) r+k(u,v,w)=0\Rightarrow $ 
\begin{equation}
r=\dfrac{k(u,v,w)}{\allowbreak 2v\left( 1-2u\right) \left( 1-w\right) }
\label{rk1}
\end{equation}

This shows that $r$ is unique and positive, by Lemma \ref{dk} and the fact
that $2v\left( 2u-1\right) \left( w-1\right) >0$ if $(u,v,w)\in Z_{1}\cup
Z_{2}\cup Z_{3}\subset X_{4}$. As noted in the proof of Lemma \ref{L1}, $%
\dfrac{s}{r}=\dfrac{\allowbreak 4(1-u)v(1-w)}{-4vw+4vwu+1}$. Now $%
4(1-u)v(1-w)-\left( -4vw+4vwu+1\right) =\allowbreak 4v-4uv-1>0$ by (\ref{5a}%
), which implies that $\dfrac{\allowbreak 4(1-u)v(1-w)}{-4vw+4vwu+1}>1$ by (%
\ref{5b}), and thus $s>r$. Using (\ref{rk1}) yields $s=\dfrac{\allowbreak
2(1-u)vk(u,v,w)}{\left( -4vw+4vwu+1\right) v\left( 1-2u\right) }$, which
also shows that $s$ is unique. Thus we have a solution $(r,s,u,v,w)$ of (\ref%
{1})--(\ref{3}) with $0<r<s$. Let $x_{1}=u-1,x_{2}=rv,$ and $x_{3}=(s-r)w+r$%
. Then (\ref{symm4}) must hold since (\ref{1})--(\ref{3}) and (\ref{symm4})
are an equivalent system of equations. If $p(x)=(x+1)x(x-r)(x-s)$, then $%
x_{1},x_{2},$ and $x_{3}$ must be the critical points of $p$ by Lemma \ref%
{cp}. Since $u=\dfrac{x_{1}-(-1)}{0-(-1)}$, $v=\dfrac{x_{2}-0}{r-0}$, and $w=%
\dfrac{x_{3}-r}{s-r}$, $(u,v,w)$ is a ratio vector of $p$.
\end{proof}

\begin{remark}
Equivalent necessary and sufficient conditions for $(u,v,w)$ to be a ratio
vector is $(u,v,w)\in Z\cap X_{4}\cap \left\{ (u,v,w):\dfrac{1}{4(1-u)}<v<%
\dfrac{1}{4(1-u)w}\right\} \cap \left\{ (u,v,w):k(u,v,w)>0\right\} $. We
preferrred to use the sets $Z_{1},Z_{2},Z_{3}$ instead since the set where $%
k(u,v,w)>0$ is not easy to determine. However, given $(u,v,w)\in Z$, one can
easily check if $k(u,v,w)>0$.
\end{remark}

\begin{remark}
One can also write the solution for $s$ as $s=\dfrac{4v(4u-1)(1-u)(1-w)}{%
d(u,v,w)}$, where $d$ is the polynomial defined in the proof of Lemma \ref%
{dk}.
\end{remark}

\textbf{Example: }Let $u=\dfrac{15}{32}\approx .468\,8$, $v=\dfrac{5}{9}%
\approx 0.555\,6$, and $w=\dfrac{156\,303-9\sqrt{10\,054\,801}}{211\,888}%
\approx .6030$. Then $\dfrac{1}{4(1-u)}\approx .470\,6<v<\dfrac{1}{4(1-u)w}%
\approx .780\,4$, and by Proposition \ref{sufficient}, $(u,v,w)$ is a ratio
vector of $p(x)=(x+1)x(x-r)(x-s)$, where $r\approx \allowbreak
5.\,\allowbreak 982\,1$ and $s\approx \allowbreak 9.\,\allowbreak 730\,5$

{\Large Further Discussion}

Let $S$ be the surface $R(u,v,w)=\allowbreak 0$. Then $S$ contains the
family of lines $u=C$, $v=\dfrac{1}{2}$, $w=1-C$, $\dfrac{1}{4(1-u)}%
=\allowbreak \dfrac{1}{4\left( 1-C\right) }<v=\dfrac{1}{2}\iff C<\dfrac{1}{2}
$, $\dfrac{1}{4(1-u)w}=\allowbreak \dfrac{1}{4\left( 1-C\right) ^{2}}>v=%
\dfrac{1}{2}\iff 1-\dfrac{1}{2}\sqrt{2}<C<1$ or $1<C<1+\dfrac{1}{2}\sqrt{2}$%
, and $k(u,v,w)=\allowbreak C\left( -1+4C-2C^{2}\right) >0\iff 0<C<1+\dfrac{1%
}{2}\sqrt{2}$ or $C<0$.

By the remark above, $\left( C,\dfrac{1}{2},1-C\right) $ is a ratio vector
if and only if $1-\dfrac{1}{2}\sqrt{2}<C<\dfrac{1}{2}$. However, we can
prove more. Since $R(u,1/2,w)=\allowbreak \left( 1-w-u\right) \left(
2wu-2w+1\right) $ and $2uw-2w+1>0$ by (\ref{7a}), $\left( u,\dfrac{1}{2}%
,w\right) $ is a ratio vector only if $u+w=1$. Thus we have proven

\begin{theorem}
$\left( u,\dfrac{1}{2},w\right) $ is a ratio vector if and only if $u=C$ and 
$w=1-C,$ where $1-\dfrac{1}{2}\sqrt{2}<C<\dfrac{1}{2}$.
\end{theorem}

Note that though S contains a family of lines, $S$ is not a ruled surface in
general. That is easy to see by looking at the second partials $\dfrac{%
\partial ^{2}R}{\partial u^{2}},$ $\dfrac{\partial ^{2}R}{\partial v^{2}},%
\dfrac{\partial ^{2}R}{\partial w^{2}}$.

\textbf{Correction: }We make a minor correction to equation (11) in (\cite{h}%
). It should read

$\left( \sigma _{1}(\sigma _{2}-\sigma _{3}+1)-\sigma _{2}+\sigma _{3}-\frac{%
1}{2}\right) r_{1}r_{3}$ $+\left( \sigma _{1}\sigma _{3}-\sigma _{3}+\frac{1%
}{2}\right) r_{1}r_{4}+\sigma _{2}(\sigma _{3}-1)r_{3}^{2}+\left( \frac{1}{2}%
-\sigma _{2}\sigma _{3}\right) r_{3}r_{4}=0$

\end{document}